\documentclass[letterpaper, 10 pt, conference]{ieeeconf}  % Comment this line out if you need a4paper

\IEEEoverridecommandlockouts                              % This command is only needed if 
\usepackage{cite}
\usepackage{amsmath,amssymb,amsfonts}
\usepackage{graphicx}
\usepackage{textcomp}
\usepackage{xcolor}
\usepackage[english]{babel}
\usepackage{tabularx}
\usepackage{algorithm}
\usepackage{algpseudocode}
\usepackage{caption}
\usepackage{subcaption}

\newtheorem{assumption}{Assumption}
\newtheorem{definition}{Definition}
\newtheorem{remark}{Remark}

\def\BibTeX{{\rm B\kern-.05em{\sc i\kern-.025em b}\kern-.08em
    T\kern-.1667em\lower.7ex\hbox{E}\kern-.125emX}}

\title{\LARGE \bf
A Physics-Based Safety Recovery Approach for Fault-Resilient Multi-Quadcopter Coordination
}

\author{Hamid Emadi, Harshvardhan Uppaluru, and  Hossein Rastgoftar
\thanks{The authors with the Aerospace and Mechanical Engineering Department at University of Arizona Emails: \{hamidemadi, huppaluru, hrastgoftar\}@email.arizona.edu}% <-this % stops a space
}

\begin{document}

\maketitle
\thispagestyle{empty}
\pagestyle{empty}

%%%%%%%%%%%%%%%%%%%%%%%%%%%%%%%%%%%%%%%%%%%%%%%%%%%%%%%%%%%%%%%%%%%%%%%%%%%%%%%%
\begin{abstract}
This paper develops a novel physics-based approach for fault-resilient multi-quadcopter coordination in the presence of abrupt quadcopter failure. Our approach consists of two main layers: (i) high-level physics-based guidance  to safely plan the desired recovery trajectory for every healthy quadcopter and (ii) low-level trajectory control design by choosing an admissible control for every  healthy quadcopter to safely recover from the anomalous situation, arisen from quadcopter failure, as quickly as possible. For the high-level trajectory planning, first, we consider healthy quadcopters as particles of an irrotational fluid flow sliding along streamline paths wrapping failed quadcopters in the shared motion space. We then  obtain the desired recovery trajectories  by maximizing the sliding speeds along the streamline paths such that the rotor angular speeds of healthy quadcopters do not exceed  certain upper bounds at all times during the safety recovery. In the low level, a feedback linearization control is designed for every healthy quadcopter such that quadcopter rotor angular speeds remain bounded and satisfy the corresponding safety constraints. Simulation results are given to illustrate the efficacy of the proposed method.

\end{abstract}

%62467
%%%%%%%%%%%%%%%%%%%%%%%%%%%%%%%%%%%%%%%%%%%%%%%%%%%%%%%%%%%%%%%%%%%%%%%%%%%%%%%%
\section{INTRODUCTION}
Unmanned aerial vehicle (UAV) was originally developed and used for military missions~\cite{peng2009design}. However, recently, applications of UAVs have been extended in different fields. For instance, Multi quadcopter systems (MQS) have been used for data acquisition from hazardous environments or agricultural farm fields, surveillance applications, urban search and rescue, wildlife monitoring and exploration \cite{argrow2005uav} \cite{tsouros2019review} \cite{witczuk2018exploring}.   One of the main notions in networked cooperative systems is fault resilient~\cite{rastgoftar2020fault}\cite{leblanc2012resilient}\cite{dibaji2017resilient}. In this work, we propose a novel physics-based approach for recovery planning of an MQS under failure of group of agents.  

\subsection{Related Work}
Multi-agent coordination is one of the main challenges in UAV-based systems. Researchers have proposed different multi-agent coordination approaches in the past. For example, authors in \cite{li2018nonlinear} proposed nonlinear consensus-based control strategies for a group of agents under different communication topologies. Another approach is containment control in which a group of followers are coordinated by a group of leaders through local communications. Authors in \cite{cao2012distributed} \cite{ji2008containment} provide distributed containment control of a group of mobile autonomous agents with multiple stationary or dynamic leaders under both fixed and switching directed network topologies. Authors in \cite{kim2008pde},\cite{krishnan2018distributed} and \cite{frihauf2010leader} propose  partial differential equations (PDE) based methods in which the position of the agents is the state of the PDE. Another coordination approach is continuum deformation proposed in \cite{rastgoftar2016continuum} \cite{rastgoftar2017continuum} \cite{rastgoftar2019physics}. This method is also based on the local communication between a group of followers a group of leaders. Graph rigidity method is proposed by \cite{wang2018distance} for the leaderless case and the
leader-follower case.  

One of the main goals in MQS is to avoid collision when an unexpected obstacle emerges in the airspace. For instance, when a quadcopter fails, the rest of quadcopters must change their path accordingly to satisfy the safety conditions. Therefore, each quadcopter must have sense and avoid (SSA) capabilities to avoid collision in case of pop up failures of other agents. Many researches have been conducted on autonomous collision avoidance of MQS. Authors in \cite{kang2017sense} propose the collision avoidance method based on estimating and predicting the agents' trajectory. A reference SAA system architecture is presented based on Boolean Decision Logics in \cite{ramasamy2015unified}. Authors in \cite{yu2015sense} provide a complete survey on SSA technologies in the sequence of fundamental functions/components of SSA in sensing techniques, decision making, path planning, and path following.  

In~\cite{rastgoftar2019physics}, authors develop a continuum deformation framework for traffic coordination management in a finite motion space. In particular, authors propose macroscopic coordination planning based on Eulerian continuum mechanics, and microscopic path planning of quadcopters considered as particles of a rigid body. This work lies in a similar vein. In this paper, we extend the work in \cite{rastgoftar2019physics} to address the scenario in which a set of failures of quadcopters are reported. We develop a physics-based approach for recovery planning, and we verify the proposed method on dynamics of a group of quadcopters.

%\begin{figure}[ht]
%    \centering
%    \includegraphics[width=0.4\textwidth]{plots/3D.png}
%    \caption{Quadcopter team configurations: blue at $t_0=0 \sec$, green at $t_0=20 \sec$, yellow at $t_0=40 \sec$ and magenta at $t_0=60 \sec$. Unsafe zone due to failure $f_1$ and $f_2$ are shown as red cylinders.}
%    \label{fig:initial conf}
%\end{figure}

\begin{figure}
     \centering
     \begin{subfigure}[b]{0.5\textwidth}
         \centering
         \includegraphics[width=\textwidth]{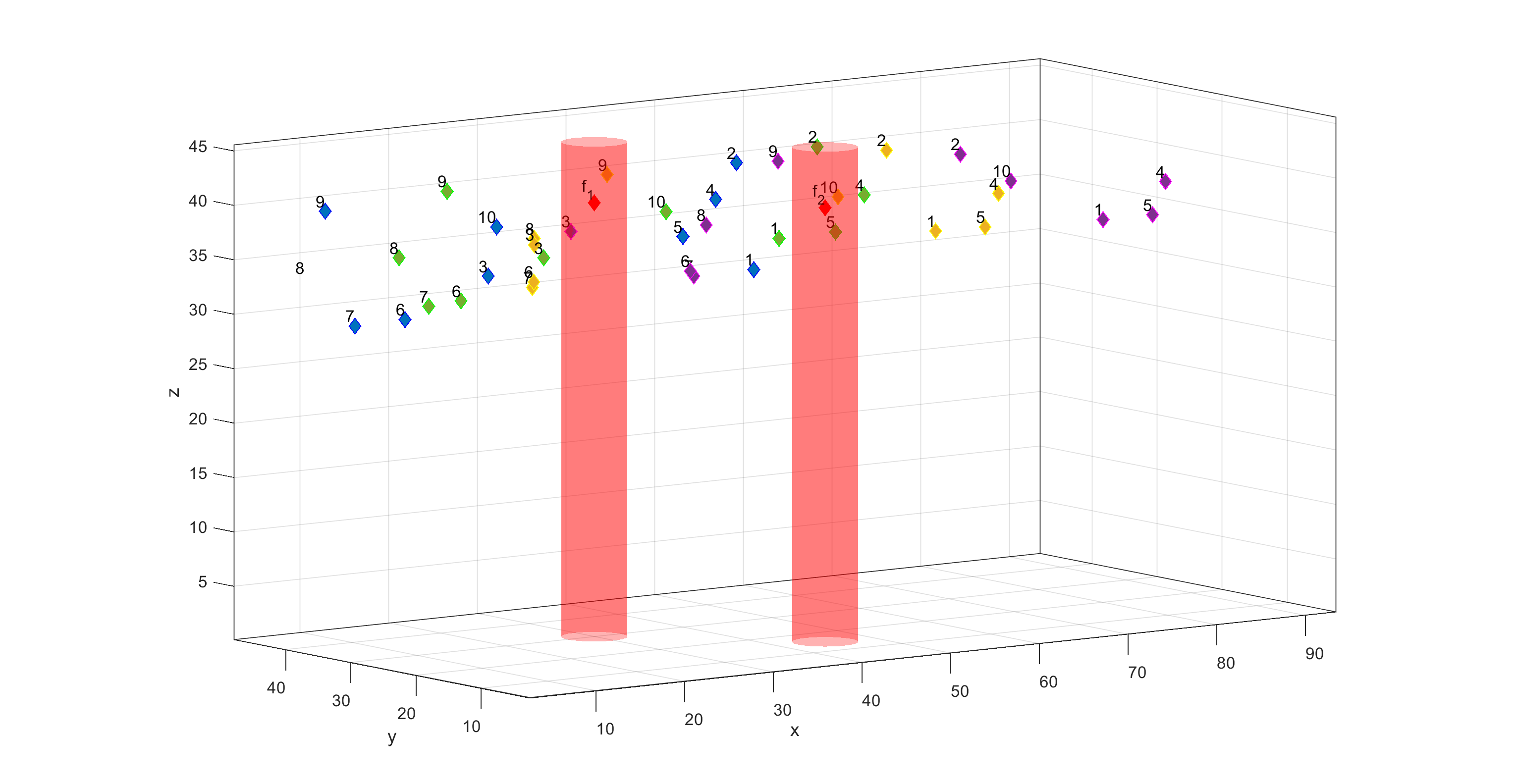}
         \vspace{-0.5cm}
         \caption{}
         \label{fig:initial conf}
     \end{subfigure}
     \hfill
     \begin{subfigure}[b]{0.5\textwidth}
         \centering
         \includegraphics[width=\textwidth]{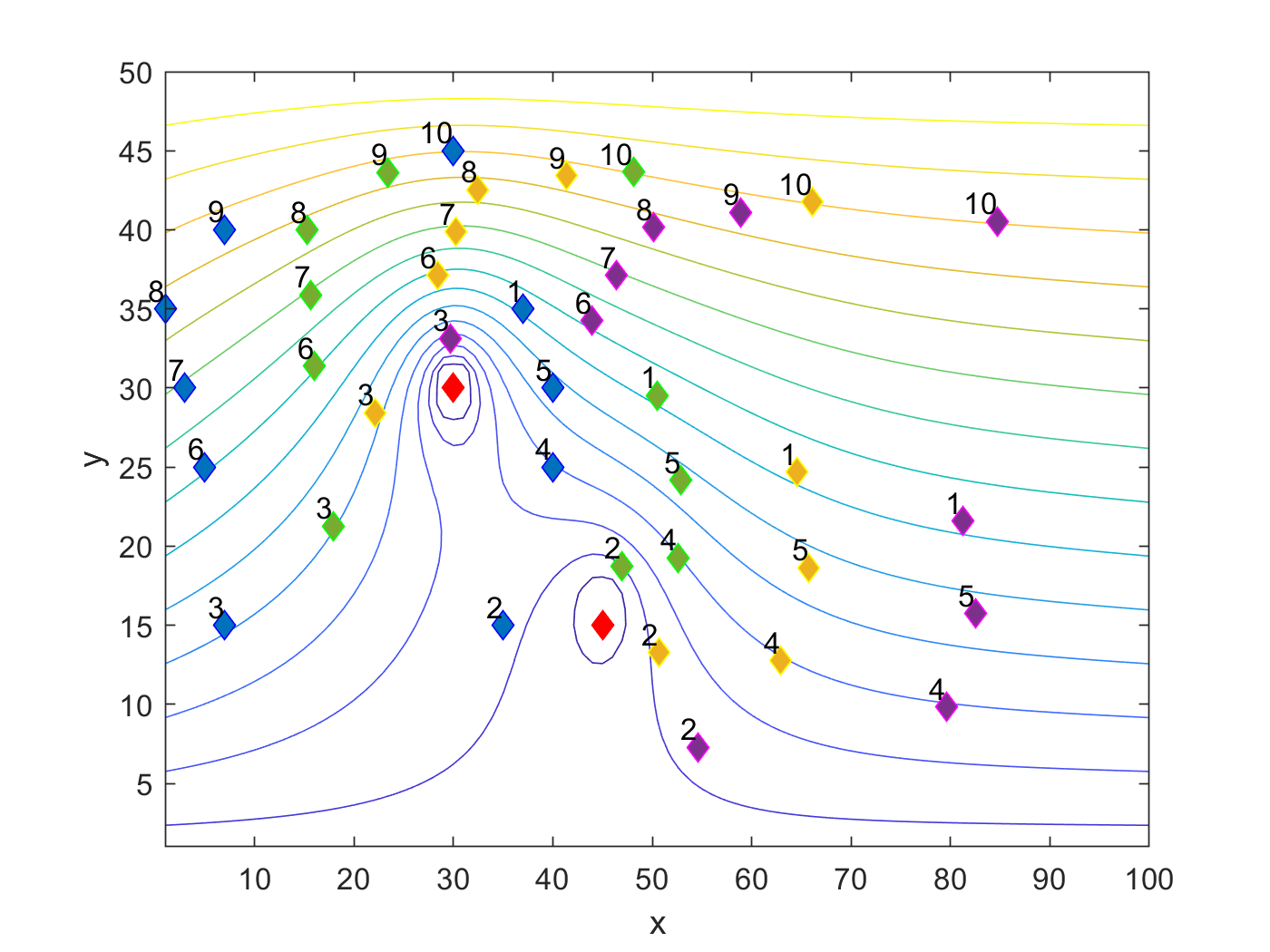}
         \vspace{-0.5cm}
         \caption{}
         \label{fig:stream_contour}
     \end{subfigure}
     
        \caption{(a) Quadcopter team configurations: blue at $t_0=0 \sec$, green at $t_0=20 \sec$, yellow at $t_0=40 \sec$ and magenta at $t_0=60 \sec$. Unsafe zone due to failure $f_1$ and $f_2$ are shown as red cylinders. (b) Stream lines in $x-y$ plane.}
        \label{fig:three graphs}
\end{figure}

\subsection{Contributions and Outline}
We propose a new physics-based approach for resilient multi-UAV coordination in the presence of UAV failure. Without loss of generality, this paper considers each UAV to be a quadcopter modeled by a $14$-th order nonlinear dynamics presented in \cite{rastgoftar2021safe}. In particular, we consider a single quadcopter team coordinating in a $3$-D motion space, and classify individual quadcopters as \textit{healthy} and \textit{failed} agents. While the healthy quadcopters can admit the desired group coordination, the failed quadcopters cannot follow the desired group coordination. To deal with this anomalous situation, we ensure safety of the healthy quadcopters and inter-agent collision avoidance by developing a two-fold safety recovery approach with planning and control layers.  For the planning of safety recovery, we treat the healthy quadcopters as particles of an ideal fluid flow field sliding along the streamline paths wrapping the failed quadcopters. For every healthy quadcopter,  the desired recovery trajectory is safely planned  by maximizing the sliding speed of the quadcopter, along the safety recovery path, such that the constraints on quadcopter rotor angular speeds are all satisfied. This safety recovery planning is complemented by designing a nonlinear recovery trajectory control for each healthy quadcopter that assures satisfaction of all safety constraints.

This paper is organized as follows: Problem Statement is discussed in Section \ref{Problem Statement}. Safety recovery planning and control are presented in Sections \ref{High-Level Planning: Recovery Trajectory Planning} and \ref{Mathematical Modeling of Quadcopters and Trajectory Tracking Control}, respectively. The results of the safety recovery simulation are presented in Section \ref{Simulation Results} and followed by Conclusion in Section \ref{Conclusion}.

\section{Problem Statement}\label{Problem Statement}
% This paper develops a framework for safety recovery of multi-quadcopter coordination in the presence of abrupt quadcopter failure. 
We consider an MQS consisting of $n_q$ quadcopters defined by set $\mathcal{I}=\left\{1,\cdots,n_q\right\}$. 
We assume that $n_f<n_q$ quadcopters identified by set $\mathcal{F}\subset \mathcal{I}$ unpredictably fail to follow the desired group coordination at reference time $t_0$ but the remaining quadcopters, defined by set $\mathcal{L}=\mathcal{I}\setminus \mathcal{F}$, can still move cooperatively and follow the desired group coordination. To safely recover from this anomalous situation, we propose to treat the healthy quadcopters as particles of an ideal fluid flow, defined by combining uniform flow in the $x-y$ plane and doublet flow. To this end, we use complex variable $\mathbf{z}=x+\mathbf{i}y$ to denote the position in the $x-y$ plane, and obtain the potential function ${\Phi}\left(x,y\right)$ and stream function ${\Psi}\left(x,y\right)$ of the ideal fluid flow field by defining 
    
\begin{eqnarray}
    \nonumber
f\left(\mathbf{z}\right)&=&{\Phi}\left(x,y\right)+\mathbf{i} {\Psi}\left(x,y\right)
    \\&=&u_{\infty}\sum_{h\in \mathcal{F}}\left(\mathbf{z}-\mathbf{z}_h+\dfrac{a_h^2}{\mathbf{z}-\mathbf{z}_h}\right),
    \label{eq: complex function}
\end{eqnarray}
over the complex plane $\mathbf{z}$, where $\mathbf{z}_h$ denotes position of the failed quadcopter $h\in \mathcal{F}$; $u_{\infty}$ and $a_h$ are constant design parameters for planning the safety recovery.

By using the ideal fluid flow model, $x$ and $y$ components of every cooperative quadcopter $i\in \mathcal{I}$ are constrained to  slide along the stream curve 
 $\Psi_i=\Psi(x_i(t),y_i(t))=\Psi_{i,0}$ at any time $t\geq t_0$, where 
 \begin{equation}\label{eq: psi initial}
     \Psi(x_i(t_0),y_i(t_0))=\Psi_{i,0},\qquad \forall i\in \mathcal{L}.
 \end{equation}
Also, every failed quadcopter is excluded from the motion space by a circular cylinder elongated in $z$ direction (see Fig. \ref{fig:initial conf}). 
\begin{remark}
If only one failed UAV exists at time $t\geq t_0$, then, the cross-section of the wrapping cylinder is a circle of radius $a_h$ centered at $\mathbf{z}_h$. Otherwise (i.e. $\left|\mathcal{F}\right|>1$), the cross section of the wrapping cylinder is not an exact circle. Note that expression \eqref{eq: complex function} specifies  a conformal mapping between the $x-y$ and $\Phi-\Psi$ planes, where $\Phi(x,y)$ and $\Psi(x,y)$ satisfy the Cauchy-Riemann and Laplace equation:
\begin{eqnarray}\label{Laplac}
    \nabla^2 \Psi = 0, \quad \nabla^2 \Phi = 0
\end{eqnarray}
\end{remark}

\begin{assumption}
We assume that healthy quadcopters move sufficiently fast or the $a_h$ is chosen sufficiently large  such that the failed quadcopters do not leave the wrapping cylinders during the the safety recovery interval. 
\end{assumption}
\begin{assumption}
We assume that the recovery trajectories of all quadcopters are planned such that the altitude remains constant. Thus, $z$ component of velocity is 0.  
\end{assumption}

By the above problem setting, the main objective of this paper is to plan the recovery trajectory for every healthy quadcopter $i\in \mathcal{L}$ so that MQS can recover safety as quickly as possible, by wrapping the failed quadcopters. Here, we assume that the rotor speeds of every quadcopter must not exceed $\omega^{max}_r$. This safety condition can be formally specified by

\begin{equation}\label{mainsafety}
    0<{\omega_r}_{i,j}(t)\leq \omega^{max}_r,\quad \forall i\in \mathcal{I},~j\in \left\{1,\cdots,4\right\},~\forall t\geq t_0
\end{equation}
where ${\omega_r}_{i,j}(t)$ is the angular speed of rotor $j\in \left\{1,\cdots,4\right\}$ of quadcopter $i\in \mathcal{L}$ at time $t\geq t_0$. $\boldsymbol{r}_i(t)$ and $\boldsymbol{r}_{i,d}(t)$ denote the actual position and desired trajectory of quadcopter $i$ at $t\geq t_0$, respectively.  
\begin{figure*}[ht]
    \centering
    \includegraphics[width=.9\textwidth]{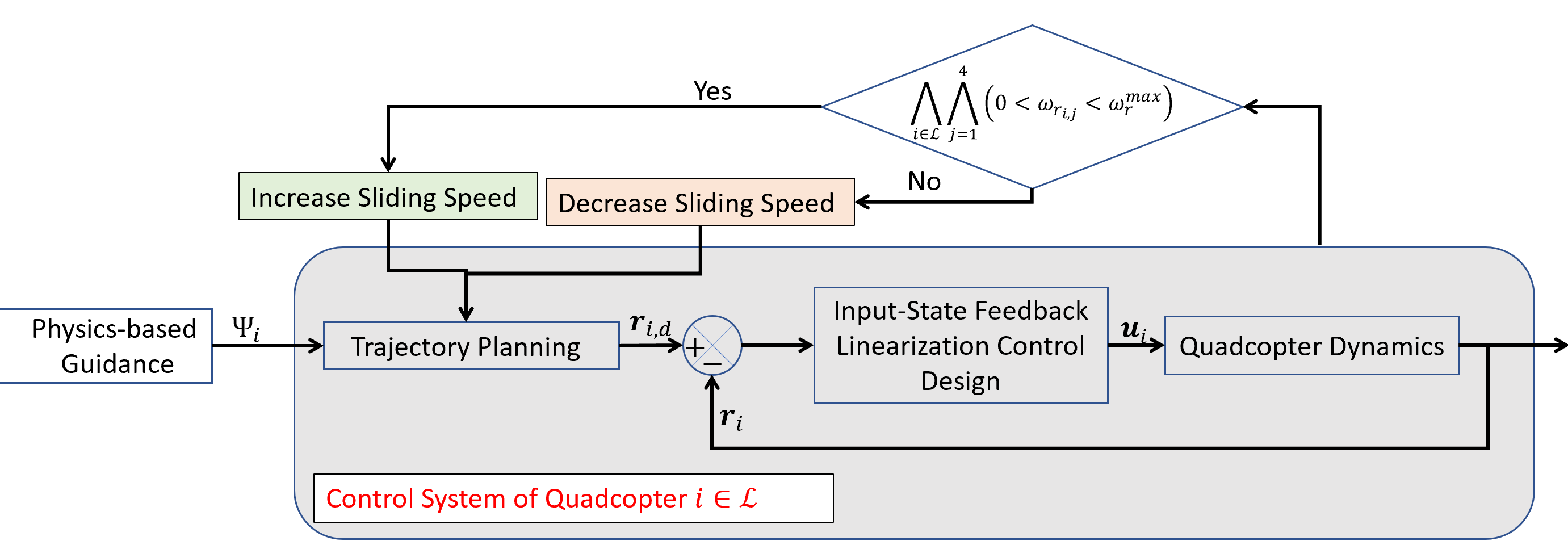}
    %\vspace{-0.5cm}
    \caption{Block diagram of MQS with the proposed method}
    \label{fig:block diagram}
\end{figure*}

We decompose this safety recovery planning into (i) high-level trajectory planning presented in Section \ref{High-Level Planning: Recovery Trajectory Planning} and (ii) low-level trajectory tracking control presented in Section  \ref{Mathematical Modeling of Quadcopters and Trajectory Tracking Control}. More specifically, Section  \ref{High-Level Planning: Recovery Trajectory Planning} obtains the safety recovery stream lines ($\Psi(x_i(t_0),y_i(t_0))=\Psi_{i,0}$)  for every healthy  quadcopter $i\in \mathcal{L}$, numerically, by using the finite difference method. This is complemented  by determining  the desired safety recovery trajectory through assignment the maximum sliding speed along the stream $\Psi_{i,0}$ ($\forall i\in \mathcal{L}$), satisfying safety condition \eqref{mainsafety}, in Section \ref{High-Level Planning: Recovery Trajectory Planning}. Section  \ref{Mathematical Modeling of Quadcopters and Trajectory Tracking Control} applies the feedback linearization control approach presented in \cite{rastgoftar2021safe} to safely track the recovery trajectory by choosing an admissible quadcopter control satisfying safety constraint \eqref{mainsafety}. Fig~\ref{fig:block diagram} shows the block diagram of MQS with the proposed approach.

\section{High-Level Planning: Recovery Trajectory Planning}\label{High-Level Planning: Recovery Trajectory Planning}
The complex function $f(\mathbf{z})$, expressed in~\eqref{eq: complex function}, provides a closed form solution for $\Phi$ and $\Psi$. However, as mentioned in Remark 1, in case of multiple failures, the area enclosed by each unsafe zone is not an exact circle, and we cannot arbitrarily shape the enclosing unsafe  area for multiple failures in the motion space. To deal with this issue,  we use the finite difference approach to determine $\Phi$ and $\Psi$ values over the motion space and arbitrarily  shape of the area enclosing the failed quadcopter.

% consequently, recovery planning deviates from considering the unsafe zones by circular cylinders. This leads to the numerical approach for computing $\Phi$ and $\Psi$. Moreover, numerical solution gives us the opportunity to deal with arbitrary geometry of unsafe zones.   

%In order to trajectory planning for a group of UAVs in case of a set of failure of agents and therefore make the airspace safe for flight of the rest of agents, we deploy the physics-based continuum deformation approach and the potential flow theory. In particular, we consider the navigable channels of UAVs as the stream lines of external flows around circular cylinders (i.e. the unsafe zone around the pop-up failure).
{
Let $\mathcal{C}$ be a set representing the projection of the airspace on the $x-y$ plane, and failures identified by set  $\mathcal{F}=\{n_1,\dots,n_f\}$  occur at $\boldsymbol{r}_{n_1},\dots,\boldsymbol{r}_{n_f}$ in $\mathcal{C}$. In the presence of abrupt quadcopter failure, quadcopters' trajectories should be modified accordingly to provide a safe maneuver in $\mathcal{C}$ and safely wrap the unsafe zones in  $\mathcal{C}$. {To this end, the unsafe zone $\mathcal{U}_i$ corresponding  to the failed quadcopter $n_i\in \mathcal{F}$ is defined by a circle with radius $a_{n_i}$  centered at $(x_{n_i},y_{n_i})$.} Then, the recovery trajectories of healthy quadcopters can be defined by the stream functions of an ideal flow around a set of circular cylinders enclosing $\mathcal{F}$. }

%Let $\Phi(\boldsymbol{r})$ denote a {\it{potential function}} which satisfies the basic laws of fluid mechanics: {\it{conservation of mass}} and {\it{momentum}}, assuming incompressible, inviscid and irrotational flow. In a potential flow, {\it{velocity vector field}}, denoted by $\boldsymbol{V}(\boldsymbol{r})$, is
%\begin{eqnarray}
%    \boldsymbol{V}(\boldsymbol{r}) = \frac{\nabla \Phi(\boldsymbol{r})}{K(\boldsymbol{r})},
%\end{eqnarray}
%where $K(\boldsymbol{r})$ is a distributed agents' velocity function that enables to treat UAVs with different speeds in the airspace. Note that $\Phi$ must satisfy the Laplace's equation (i.e conservation of mass $\nabla.\boldsymbol{V}$=0):
%\begin{eqnarray}\label{Laplac Phi}
%    \nabla^2 \Phi = 0.
%\end{eqnarray}
% From~\eqref{Cauchy Reimann 1} and~\eqref{Cauchy Reimann 2} , in two dimensions, stream function must satisfy the Laplace's equation:
%\begin{eqnarray}\label{Laplac Psi}
%    \nabla^2 \Psi = 0.
%\end{eqnarray}

%In this work, we assume that each agent's objective is to follow an admissible trajectory in order to avoid collision to the failed agents. Although, in general, $\Phi$ and $\Psi$ are time varying functions, we assume that the unsafe space, due to the failed UAVs, is a time-invariant space, and consequently, we treat with $\Phi(\boldsymbol{r})$ and $\Psi(\boldsymbol{r})$ as steady state potential and stream functions, respectively.    

\begin{figure}[ht]
    \centering
    \includegraphics[width=0.4\textwidth]{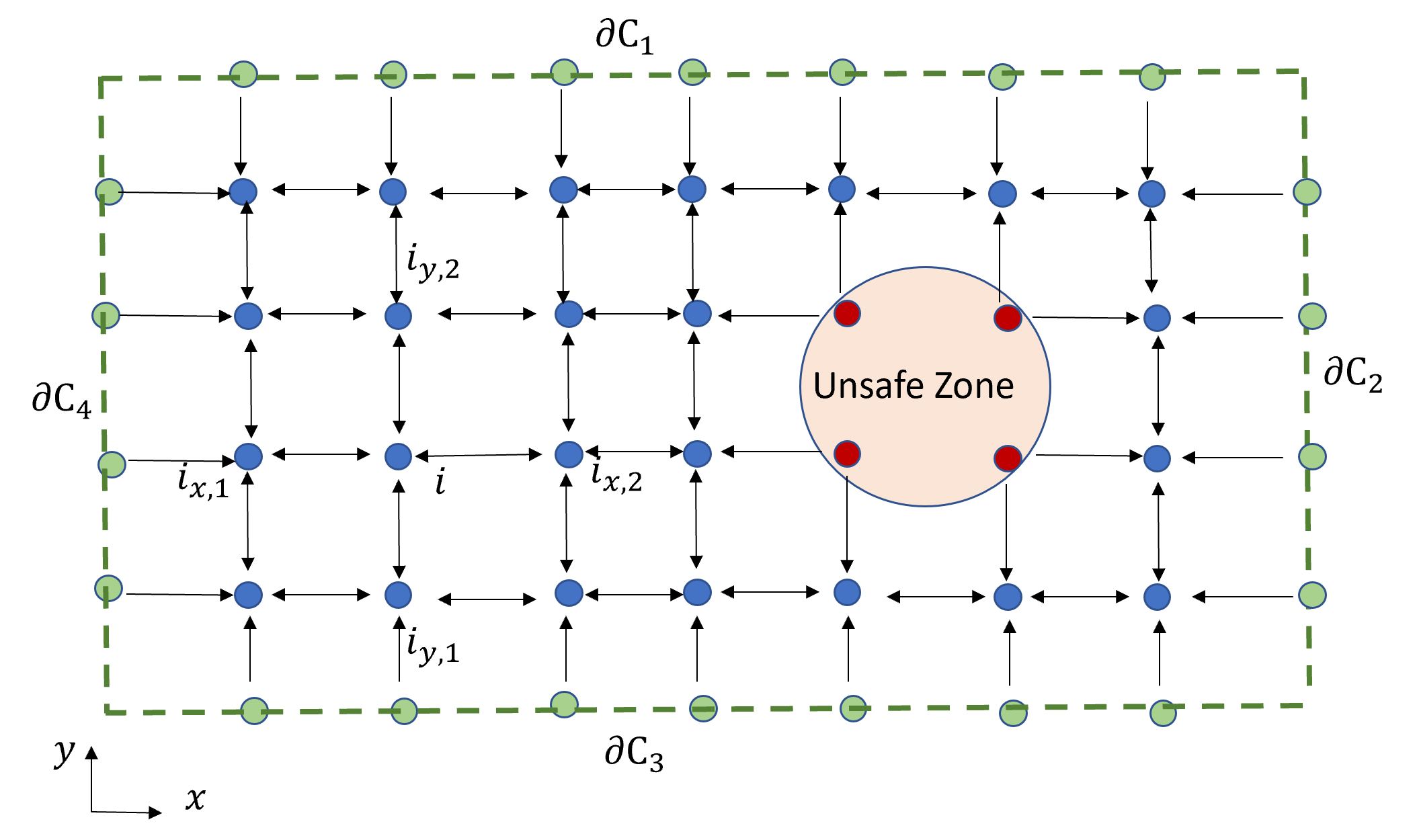}
    \caption{Directed graph $G$ resulted from discretizing $x-y$ plane }
    \vspace{-0.5cm}
    \label{fig:grid}
\end{figure}

Without loss of generality, we assume that $\mathcal{C}$ is a rectangular environment lies in the $x-y$ plane and
% along the direction of motion of quadcopters, and $i^{\text{th}}$ quadcopter aims to move from $\boldsymbol{r}_{i,0} = [x_{i,0},y_{i,0},z_{i,0}]^T$ to $\boldsymbol{r}_{i,f} = [x_{i,f},y_{i,0},z_{i,0}]^T$. Note that each quadcopter stays in the $x-y$ plane during the mission. In order to compute $\Phi(x,y),\Psi(x,y)$ over $\mathcal{C}$, \eqref{Laplac} needs to be solved with proper boundary conditions. In the following, we
use the finite difference method to compute $\Phi,\Psi$ over $\mathcal{C}$. The idea of finite-difference-method is to discretize the governing  PDE and the environment by replacing the partial derivatives with their approximations. We uniformly discretize $\mathcal{C}$ into small regions with increments in the $x,y$ directions given as $\Delta x,\Delta y$, respectively. Discretizing $\mathcal{C}$ in $x-y$ plane results in the directed graph $G(\mathcal{V},\mathcal{E})$ in which, each node is connected to the adjacent nodes in $x$ and $y$ direction (Fig~\ref{fig:grid}). Node set and edge set are defined as $\mathcal{V} =\{1,\dots,m\}$ and $\mathcal{E} \subseteq \mathcal{V} \times \mathcal{V}$, respectively. $\mathcal{E}$ is a set of pairs $(i,j)$ connecting nodes $i,j\in \mathcal{V}$.  

Without loss of generality, suppose that nodes are labeled such that the boundary nodes, interior nodes over the safe zone and interior nodes over the unsafe zone are labeled as $\mathcal{V}_b=\{1,\dots,{m_b}\}$,  $\mathcal{V}_f = \{{m_b+1},\dots,{m_b+m_f}\}$ and $\mathcal{V}_c=\{{m_b+m_f+1},\dots,{m}\}$, respectively. Let $\partial\mathcal{C}_1,\partial\mathcal{C}_3$ and and $\partial\mathcal{C}_2,\partial\mathcal{C}_4$ denote the boundaries of rectangular area $\mathcal{C}$ in $x$ and $y$ directions, respectively (see Fig~\ref{fig:grid}). We plan the safety recovery  trajectories such that the average bulk motion of the healthy  MQS is from left to right along the positive  $x$ direction. To fulfill this requirement, we choose the boundary conditions of $\Psi$ as follows:
\begin{eqnarray} \label{boundary conditions}
    \Psi(j) =  
    \left\{\begin{matrix}
       Ky_j & j\in \partial\mathcal{C}_1 \bigcup \partial\mathcal{C}_2 \bigcup\partial\mathcal{C}_3\bigcup \partial\mathcal{C}_4 \\ 
        0 &  j\in \mathcal{V}_f
    \end{matrix}\right.
\end{eqnarray}
where $y_j$ is the $y$ component of position of node $j$, and $K$ is a positive constant number. From the above expression, $\Psi$ is constant over $\partial\mathcal{C}_1 , \partial\mathcal{C}_2$. Hence, $\partial\mathcal{C}_1 , \partial\mathcal{C}_2$ are stream lines.

By substituting the approximated derivatives from the Taylor series to~\eqref{Laplac}, stream value function $\Psi_i$ at node $i\in \mathcal{V}_c$ satisfies the following equation:
\begin{eqnarray}\label{discrete_laplace_phi}
    \frac{\Psi_{i_{x,1}}-2\Psi_i+\Psi_{i_{x,2}}}{{\Delta x}^2} + \frac{\Psi_{i_{y,1}}-2\Psi_i+\Psi_{i_{y,2}}}{{\Delta y}^2} = 0,
\end{eqnarray}
where $\Psi_{i_{x,1}}$ and $\Psi_{i_{x,2}}$ are potential values at neighbor nodes in $x$ direction. Similarly, $\Psi_{i_{y,1}}$ and $\Psi_{i_{y,2}}$ are the potential values at neighbor nodes in $y$ direction.

Let $\boldsymbol{\Psi} = \begin{bmatrix}
\Psi_1 &\dots&\Psi_m
\end{bmatrix}^T$ represent the nodal vector of the potential function. \eqref{discrete_laplace_phi} can be written in the compact form of 
\begin{eqnarray}\label{eq: L * phi=0}
    L\boldsymbol{\Psi} = \boldsymbol{0}.
\end{eqnarray}
where $L\in \mathbb{R}^{m\times m}$ is the Laplacian matrix of the network. Entries of $L$ are defined as
%\begin{eqnarray}
%    L_{ij}=\left\{\begin{matrix}
%     -1&  i=j,i\in\mathcal{V}_c\bigcup\mathcal{V}_f& \\ 
%     0.25&  i\neq j,(i,j)\in \mathcal{E}, i\in\mathcal{V}_c, j\in\mathcal{V}_c& \\ 
%     0.25&  i\neq j,(i,j)\in \mathcal{E}, i\in\mathcal{V}_c, j\in\mathcal{V}_f& \\ 
%     0.25&  i\neq j,(i,j)\in \mathcal{E}, i\in\mathcal{V}_c, j\in\mathcal{V}_b& \\
%     0.25&  i\neq j,(i,j)\in \mathcal{E}, i\in\mathcal{V}_f, j\in\mathcal{V}_c& \\
%     0&  \text{otherwise}& 
%\end{matrix}\right.
%\end{eqnarray}
\begin{eqnarray}
    L_{ij}=\left\{\begin{matrix}
     \text{deg}(i)&  i=j  \\ 
     -1&  i\neq j,(i,j)\in \mathcal{E}\\
     0&  \text{otherwise}
\end{matrix}\right.
\end{eqnarray}
where $\text{deg}(i)$ is the in-degree of node $i$. According to~\cite{veerman2020primer} the multiplicity of the eigenvalue 0 of $L$ equals to the number of maximal reachable vertex sets. In other words, multiplicity of zero eigenvalues is the number of trees needed to cover $G$. Therefore, matrix $L$ has $m_b+m_f$ eigenvalues equal to 0. Hence, rank of $L$ is $m-m_b+m_f$, and \eqref{eq: L * phi=0} can be solved for unknown values of $\boldsymbol{\Psi}$ corresponding to the interior nodes.

By obtaining  $\mathbf{\Psi}$ over $\mathcal{C}$, recovery path of healthy quadcopter $i\in \mathcal{L}$ is an stream line $\Psi_i$ defined by \eqref{eq: psi initial}.
Note that the stream line $\Psi_i$ is tangent to the desired velocity of quadcopter  $i\in \mathcal{L}$. By provoking the Cauchy-Riemann Theorem, the desired velocity of quadcopter $i\in\mathcal{L}$ is given by
\begin{eqnarray}\label{eq: rdot}
    \Dot{\boldsymbol{r}}_{i,d} = v_i\left(\frac{\partial \Psi}{\partial y} \hat{\boldsymbol{i}} - \frac{\partial \Psi}{\partial x} \hat{\boldsymbol{j}}\right) + 0\hat{\boldsymbol{k}},\qquad \forall i\in \mathcal{L},
\end{eqnarray}
where $v_i$ is the sliding speed of quadcopter $i\in \mathcal{L}$. Without loss of generality, we assume that all quadcopters move with the same sliding speed $v$ during the safety recovery. Therefore,
\begin{equation}
    v_i(t)=v(t),\qquad \forall i\in \mathcal{L},~\forall t\geq t_0.
\end{equation}
To recover safety as quickly as possible, we maximize $v$ such that the safety conditions presented in \eqref{mainsafety} are all satisfied. To this end, the maximum sliding $v^*$ is assigned by bi-section method as shown in Fig. \ref{fig:block diagram}.

Consequently, by integrating from~\eqref{eq: rdot}, we can update the desired trajectories for all agents in case of existence of failure(s) in $\mathcal{C}$.

%---------------------------------------------------------------------
%Note that if a complex function $f(x+iy) = \Phi + i\Psi$ is holomorphic, then u and v have first partial derivatives with respect to x and y, and satisfy the Cauchy–Riemann equations~\eqref{Cauchy Reimann 1},\eqref{Cauchy Reimann 2} and Laplace's equations~\eqref{Laplac Phi} and~\eqref{Laplac Psi}.
%Each unsafe failure zone can be excluded from the airspace by combining the uniform flow and doublet [REFERENCE]. In this work, we use the following complex function to compute $\Phi,\Psi$: 
%\begin{eqnarray}
%    f = U_{\infty} \sum_{i=1}^{f_n}{\left[(z-z_{f_i})+\frac{a_{f_i}^2}{(z-z_{f_i})}\right]}
%\end{eqnarray}
%where $U_{\infty}$ is the flow speed far from the failures. From the above function, $\Phi$ and $\Psi$ are computed as follows:
%\begin{eqnarray}
%    \Phi = U_{\infty} \sum_{i=1}^{f_n}{\left[(x-x_{f_i})+\frac{a_{f_i}^2(x-x_{f_i})}{(x-x_{f_i})^2+(y-y_{f_i})^2}\right]}
%\end{eqnarray}

%\begin{eqnarray}\label{Psi}
%    \Psi = U_{\infty} \sum_{i=1}^{f_n}{\left[(y-y_{f_i})-\frac{a_{f_i}^2(y-y_{f_i})}{(x-x_{f_i})^2+(y-y_{f_i})^2}\right]}
%\end{eqnarray}

% \section{Mid-Level Planning}

\section{Mathematical Modeling of Quadcopters and Trajectory Tracking Control}
\label{Mathematical Modeling of Quadcopters and Trajectory Tracking Control}
\subsection{Equations of motion}
In this work, we consider the following assumptions in mathematical modeling of quadcopter motions.

\begin{assumption}
Quadcopter is a symmetrical rigid body with respect to the axes of body-fixed frame.
\end{assumption}
\begin{assumption}
Aerodynamic loads are neglected due to low speed assumption for quadcopters.
\end{assumption}

Let $\hat{\boldsymbol{i}},\hat{\boldsymbol{j}},\hat{\boldsymbol{k}}$ be the base unit vectors of inertial coordinate system, and $\hat{\boldsymbol{i}}_b,\hat{\boldsymbol{j}}_b,\hat{\boldsymbol{k}}_b$ be the base unit vectors of a body-fixed coordinate system whose origin is at the center of mass of the quadcopter. In this section, for convenience, we omitted $i$ subscript of $i^{\text{th}}$ quadcopter in the governing equations. The attitude of the quadcoper is defined by three Euler angles $\phi,\theta$ and $\psi$ as roll angle, pitch angle and yaw angle, respectively. In this work, we use 3-2-1 standard Euler angles to determine orientation of the quadcopter. Therefore, the rotation matrix between fixed-body frame and the inertial frame can be written as
\begin{eqnarray}
    \boldsymbol{R}(\phi,\theta,\psi) = \boldsymbol{R}(\phi,0,0)\boldsymbol{R}(0,\theta,0)\boldsymbol{R}(0,0,\psi)\\
    =\begin{bmatrix}
     c_{\theta}c_{\psi}&  s_{\theta}c_{\psi}s_{\phi}-s_{\psi}c_{\phi}& s_{\theta}c_{\psi}c_{\phi}+s_{\psi}s_{\phi}\\ 
     c_{\theta}s_{\psi}&  s_{\theta}s_{\psi}s_{\phi}+c_{\psi}c_{\phi}& s_{\theta}s_{\psi}c_{\phi}-c_{\psi}s_{\phi}\\ 
    -s_{\theta}&  c_{\theta}s_{\phi}& c_{\theta}c_{\phi}
    \end{bmatrix}
\end{eqnarray}
where $s(.)=\sin(.),c(.)=\cos(.)$. Let $\boldsymbol{r}=\begin{bmatrix} x&y&z \end{bmatrix}^T$ denote the position of the center of mass of the quadcopter in inertial frame, and $\boldsymbol{\omega}=\begin{bmatrix}
\omega_x&\omega_y&\omega_z
\end{bmatrix}^T$  denote the angular velocity of the quadcopter represented in the fixed-body frame.

Using the Newton-Euler formulas, equations of motion of a quadcopter can be written in the following form:
\begin{eqnarray}\label{translation_eq}
    &\Ddot{\boldsymbol{r}} = -g\hat{\boldsymbol{k}}+\frac{p}{m}\hat{\boldsymbol{k}}_b,\\\label{rotation_eq}
    &\Dot{\boldsymbol{\omega}} = -\boldsymbol{J}^{-1}\left[\boldsymbol{\omega}\times(\boldsymbol{J}\boldsymbol{\omega})\right]+\boldsymbol{J}^{-1}\boldsymbol{\tau},
\end{eqnarray}
where $m,\boldsymbol{J}$ denote, respectively, mass and mass moment of inertia of the quadcopter. $g$ is the gravity acceleration and $p$ is the thrust force generated by the four rotors. Relation between the thrust force $p$ and angular speed of the rotors, denoted by $\omega_{r_i}$, can be written as
\begin{eqnarray}
    p = \sum_{i=1}^{4}{{f_r}_i} = b\sum_{i=1}^{4}{{{\omega^2_{r_i}}}},
\end{eqnarray}
where $b$ is the aerodynamic force constant ($b$ is a function of the density of air, the shape of the blades, the number of the blades, the chord length of the blades, the pitch angle of the blade airfoil and the drag constant), and ${f_r}_i$ is the thrust force of $i^{\text{th}}$ rotor. In~\eqref{translation_eq}, $\boldsymbol{\tau}=\begin{bmatrix} \tau_\phi&\tau_\theta&\tau_\psi
\end{bmatrix}^T$ is the control torques generated by four rotors. Relation between the $\boldsymbol{\tau}$ and angular speed of the rotors can be written in the following form
\begin{eqnarray}
    \boldsymbol{\tau}=\begin{bmatrix}
     \tau_\phi \\ 
    \tau_\theta \\ 
    \tau_\psi
    \end{bmatrix}=
    \begin{bmatrix}
     bl({\omega^2_{r_4}}-{\omega^2_{r_2}}) \\ 
    bl({\omega^2_{r_3}}-{\omega^2_{r_1}}) \\ 
    k({\omega^2_{r_2}}+{\omega^2_{r_4}}-{\omega^2_{r_1}}-{\omega^2_{r_3}})
    \end{bmatrix},    
\end{eqnarray}
where $l$ is the distance of each rotor from center of the quadcopter, and $k$ is a positive constant corresponding to the aerodynamic torques. By concatenating $p$ and $\boldsymbol{\tau}$ as input vector to the system, we can write
\begin{eqnarray} \label{eq: angular rotor speed}
    \boldsymbol{u}=\begin{bmatrix}
    p\\
     \tau_\phi \\ 
    \tau_\theta \\ 
    \tau_\psi
    \end{bmatrix}=\begin{bmatrix}
    b & b & b & b\\
    0 & -bl & 0 & bl\\ 
    -bl & 0 & bl & 0 \\ 
    -k & k & -k & k
    \end{bmatrix}\begin{bmatrix}
    {\omega^2_{r_1}}\\
    {\omega^2_{r_2}} \\ 
    {\omega^2_{r_3}} \\ 
    {\omega^2_{r_4}}
    \end{bmatrix}.
\end{eqnarray}

By defining state vector $\boldsymbol{x}=\begin{bmatrix}\boldsymbol{r}^T&\Dot{\boldsymbol{r}}^T&\phi&\theta&\psi&\boldsymbol{\omega}^T
\end{bmatrix}^T$ and input vector $\boldsymbol{u}=\begin{bmatrix}
p& \tau_{\phi}&\tau_{\theta}&\tau_{\psi}
\end{bmatrix}^T$, \eqref{translation_eq},\eqref{rotation_eq} can be written in the state space non-linear form of 
\begin{eqnarray} \label{non-linear state space}
     \left\{\begin{matrix}
        \Dot{\boldsymbol{x}}= \boldsymbol{f}(\boldsymbol{x})+\boldsymbol{g}(\boldsymbol{x})\boldsymbol{u}\\ 
        \boldsymbol{r}=\boldsymbol{C}\boldsymbol{x}
    \end{matrix}\right.   
\end{eqnarray}
where, $\boldsymbol{f}(\boldsymbol{x})$ and $\boldsymbol{g}(\boldsymbol{x})$ are defined as 
\begin{eqnarray}
     \boldsymbol{f}(\boldsymbol{x})=\begin{bmatrix}
    \boldsymbol{v}\\
     -g\hat{\boldsymbol{k}} \\ 
    \boldsymbol{\Gamma^{-1}}\boldsymbol{\omega} \\ 
    -\boldsymbol{J}^{-1}\left[\boldsymbol{\omega}\times(\boldsymbol{J}\boldsymbol{\omega})\right]
    \end{bmatrix},
\end{eqnarray}

\begin{eqnarray}
    \boldsymbol{g}(\boldsymbol{x})=\begin{bmatrix}
    \boldsymbol{0}_{3\times1} & \boldsymbol{0}_{3\times3}\\
    \frac{\hat{\boldsymbol{k}}_b}{m} & \boldsymbol{0}_{3\times3}\\ 
    \boldsymbol{0}_{3\times1} & \boldsymbol{0}_{3\times3}\\ 
    \boldsymbol{0}_{3\times1} & \boldsymbol{J}^{-1}
    \end{bmatrix}
\end{eqnarray}
and $\boldsymbol{C} = [\boldsymbol{I}_{3\times3}, \boldsymbol{0}_{3\times9}]$. $\boldsymbol{v}$ is the velocity vector of the quadcopter, and $\boldsymbol{\Gamma}$ is the matrix which relates Euler angular velocity to the angular velocity of the quadcopter. $\boldsymbol{0}_{i\times j}$ is a $i\times j $ zero matrix. In order to find $\boldsymbol{\Gamma}$, we can represent $\boldsymbol{\omega}$ in the following form
\begin{eqnarray}\label{angular velocity}
    \boldsymbol{\omega} = \Dot{\psi} \hat{\boldsymbol{k}}_1+\Dot{\theta} \hat{\boldsymbol{j}}_2 + \Dot{\phi} \hat{\boldsymbol{i}}_b,
\end{eqnarray}
where $\hat{\boldsymbol{j}}_2=\boldsymbol{R}(\phi,0,0)\hat{\boldsymbol{j}}_b$ and $\hat{\boldsymbol{k}}_1=\boldsymbol{R}(\phi,\theta,0)\hat{\boldsymbol{k}}_b$. Consequently, 
\begin{eqnarray}
    \boldsymbol{\Gamma} =\begin{bmatrix}
    1 & 0 &-s_{\theta}\\
    0 & c_{\phi} &c_{\theta}s_{\phi}\\ 
    0 & -s_{\phi} &c_{\phi}c_{\theta}
    \end{bmatrix} .
\end{eqnarray}
From~\eqref{angular velocity}, the angular acceleration $\Dot{\boldsymbol{\omega}}$ can be formulated in the following way:
\begin{eqnarray}\label{angular velocity 1}
    \Dot{\boldsymbol{\omega}} = \Tilde{\boldsymbol{B}}_1 \begin{bmatrix}
    \ddot{\phi}& \ddot{\theta}&\ddot{\psi}
    \end{bmatrix}^T +\Tilde{\boldsymbol{B}}_2
\end{eqnarray}
where $\Tilde{\boldsymbol{B}}_1=\begin{bmatrix}
\hat{\boldsymbol{i}}_b& \hat{\boldsymbol{j}}_2& \hat{\boldsymbol{k}}_1
\end{bmatrix}$ and
\begin{eqnarray} \label{tilde_B2}
    \Tilde{\boldsymbol{B}}_2=\Dot{\theta}\Dot{\psi}(\hat{\boldsymbol{k}}_1 \times \hat{\boldsymbol{j}}_2)+ \Dot{\phi}(\Dot{\psi}\hat{\boldsymbol{k}}_1+\Dot{\theta}\hat{\boldsymbol{j}}_2 )\times \hat{\boldsymbol{i}}_b
\end{eqnarray}
On the other hand, from~\eqref{non-linear state space},
\begin{eqnarray} \label{angular velocity 2}
    \Dot{\boldsymbol{\omega}} = \boldsymbol{J}^{-1}\left( -\boldsymbol{\omega}\times(\boldsymbol{J}\boldsymbol{\omega}) + \begin{bmatrix}
    u_{2}& u_3 &u_4
    \end{bmatrix}^T \right).
\end{eqnarray}
From~\eqref{angular velocity 1} and \eqref{angular velocity 2}
\begin{eqnarray}\label{equation: u to ddot}
    \begin{bmatrix}
    u_2\\ u_3\\ u_4
    \end{bmatrix} = \boldsymbol{J}\Tilde{\boldsymbol{B}}_1\begin{bmatrix}
     \ddot{\phi}\\ \ddot{\theta}\\\ddot{\psi} 
    \end{bmatrix} + \boldsymbol{J}\Tilde{\boldsymbol{B}}_2+\boldsymbol{\omega}\times(\boldsymbol{J}\boldsymbol{\omega})
\end{eqnarray}

\subsection{Recovery control}
In this subsection, we provide the input control for the non-linear state space system~\eqref{non-linear state space} to track the desired trajectory $\boldsymbol{r_d}$ obtained from section \ref{High-Level Planning: Recovery Trajectory Planning}. Since we consider low speed quadcopters, agents have enough time to update their path in case of failures. Moreover, we suppose $\boldsymbol{r}_d$ is a smooth function for all  $t\geq t_0$ (i.e. $\boldsymbol{r}_d$ has derivatives of all orders). 

In this work, we use the input-output feedback linearization approach\cite{slotine1991applied} to design the input control for a quadcopter to track the desired trajectory~\cite{rastgoftar2021safe}. We use the Lie derivative notation which is defined in the following.
\begin{definition}
Let $h:\mathbb{R}^n\rightarrow  \mathbb{R}$ be a smooth scalar function, and $\boldsymbol{f}:\mathbb{R}^n\rightarrow  \mathbb{R}^n$ be a smooth vector field on $\mathbb{R}^n$. Lie derivative of $h$  with respect to $\boldsymbol{f}$ is a scalar function defined by $L_{\boldsymbol{f}} h=\nabla h\boldsymbol{f}$.
\end{definition}

Concept of input-output linearization is based on differentiating the output until the input appears in the derivative expression. Since $u_2,u_3$ and $u_4$ do not appear in the derivative of outputs, we use the technique, called dynamic extension, in which we redefine the input vector $\boldsymbol{u}$ as the derivative of some of the original system inputs. In particular, we define $\Tilde{\boldsymbol{x}}=\begin{bmatrix}\boldsymbol{x}^T&p&\Dot{p}
\end{bmatrix}^T$ and $\Tilde{\boldsymbol{u}}=\begin{bmatrix}
u_p&\tau_{\phi}&\tau_{\theta}&\tau_{\psi}
\end{bmatrix}^T$. Therefore, extended dynamics of the quadcopter can be expressed in the following form \cite{rastgoftar2021safe}:
\begin{eqnarray} \label{extended non-linear state space}
     \left\{\begin{matrix}
        \Dot{\Tilde{\boldsymbol{x}}}= \Tilde{\boldsymbol{f}}(\Tilde{\boldsymbol{x}})+\Tilde{\boldsymbol{g}}(\Tilde{\boldsymbol{x}})\Tilde{\boldsymbol{u}}\\ 
        \boldsymbol{r}=\Tilde{\boldsymbol{C}}\Tilde{\boldsymbol{x}}
    \end{matrix}\right.   
\end{eqnarray}
where, $\Tilde{\boldsymbol{f}}(\Tilde{\boldsymbol{x}})$ and $\Tilde{\boldsymbol{g}}(\Tilde{\boldsymbol{x}})$ are defined as 
\begin{eqnarray}
     \Tilde{\boldsymbol{f}}(\Tilde{\boldsymbol{x}})=\begin{bmatrix}
    \boldsymbol{f}(\boldsymbol{x})\\
     \Dot{p} \\ 
    0\end{bmatrix}+
    \begin{bmatrix}
    \boldsymbol{0}_{3\times1}\\
     \frac{p}{m}\hat{\boldsymbol{k}}_b \\ 
    \boldsymbol{0}_{8\times1}\end{bmatrix},
\end{eqnarray}

\begin{eqnarray}
    \Tilde{\boldsymbol{g}}(\Tilde{\boldsymbol{x}})=\begin{bmatrix}
    \boldsymbol{0}_{9\times1} & \boldsymbol{0}_{9\times3}\\
    \boldsymbol{0}_{3\times1} & \boldsymbol{J}^{-1}\\ 
    0 & \boldsymbol{0}_{1\times3}\\ 
    1 & \boldsymbol{0}_{1\times3}
    \end{bmatrix}.
\end{eqnarray}
Let $\Tilde{\boldsymbol{g}}_i(\Tilde{\boldsymbol{x}})$ denote the $i^\text{th}$ column of matrix $\Tilde{\boldsymbol{g}}(\Tilde{\boldsymbol{x}})$ and $\Tilde{\boldsymbol{u}}=\begin{bmatrix}\Tilde{u}_1\dots\Tilde{u}_4
\end{bmatrix}^T$ where $\Tilde{u}_1,\dots,\Tilde{u}_4$ corresponds to $u_p,\tau_{\phi},\tau_{\theta},\tau_{\psi}$, respectively. We consider the position of the quadcopter as the output of the system (i.e. $x,y,z$). Inputs appear in the fourth order derivative of the outputs. particularly, for $q\in \{x,y,z\}$ 
\begin{eqnarray}
    \ddddot{q} = L^4_{\Tilde{\boldsymbol{f}}}q + \sum_{i=1}^{4}{L_{\Tilde{\boldsymbol{g}_i}}{L_{\Tilde{\boldsymbol{f}}}}^3q}\Tilde{u}_i
\end{eqnarray}
where $L_{\Tilde{\boldsymbol{g}_i}}{L_{\Tilde{\boldsymbol{f}}}}^3q \neq 0$ for $i=1,\dots,4$. By choosing the state transformation $\mathcal{T}(\Tilde{\boldsymbol{x}})=\begin{bmatrix}
\boldsymbol{z}&\boldsymbol{\zeta}
\end{bmatrix}^T$, \eqref{extended non-linear state space} can be converted to the following internal and external dynamics:
\begin{eqnarray}\label{internal linear dynamics}
    \Dot{\boldsymbol{\zeta}} = \begin{bmatrix}
    0 & 0\\
    0 & 1 \end{bmatrix} \boldsymbol{\zeta} +
    \begin{bmatrix}
    0\\
    1 \end{bmatrix} u_{\psi}
\end{eqnarray}
\begin{eqnarray}\label{external linear dynamics}
    \Dot{\boldsymbol{z}} = \boldsymbol{A}\boldsymbol{z}+\boldsymbol{B}\boldsymbol{s}
\end{eqnarray}
where $\boldsymbol{z} = \begin{bmatrix}
\boldsymbol{r}^T& \Dot{\boldsymbol{r}}^T& \Ddot{\boldsymbol{r}}^T& \dddot{\boldsymbol{r}}^T, 
\end{bmatrix}^T$, and  $\boldsymbol{\zeta} = \begin{bmatrix}
\psi&\Dot{\psi}
\end{bmatrix}^T$
\begin{eqnarray}
    \boldsymbol{A}=\begin{bmatrix}
    \boldsymbol{0}_{9\times 3} & \boldsymbol{I}_{9}\\
    \boldsymbol{0}_{3\times 3} & \boldsymbol{0}_{3\times 9}\end{bmatrix}, \boldsymbol{B}=\begin{bmatrix}
    \boldsymbol{0}_{9\times 3} \\
    \boldsymbol{I}_{3} \end{bmatrix}
\end{eqnarray}
where $\boldsymbol{I}_{j}$ is a $j\times j$ identity matrix. 

Next, we can figure out the Control inputs $\boldsymbol{s}$ and $u_\psi$, such that the linear systems \eqref{internal linear dynamics} and \eqref{external linear dynamics} track the desired trajectory $\boldsymbol{r}_d$. By choosing
\begin{eqnarray}
    u_{\psi} = -K_1\Dot{\psi}-K_2\psi
\end{eqnarray}
where $K_1>0,K_2>0$. Thus, the internal dynamics~\eqref{internal linear dynamics} asymptotically converges to $\boldsymbol{0}$. Moreover, we choose 
\begin{eqnarray}
    \boldsymbol{s}=K_3\left( \dddot{\boldsymbol{r}}_d-\dddot{\boldsymbol{r}}\right)+ K_4\left( \ddot{\boldsymbol{r}}_d-\ddot{\boldsymbol{r}}\right)+ \\ \nonumber
    K_5\left( \dot{\boldsymbol{r}}_d-\dot{\boldsymbol{r}}\right) + K_6\left( {\boldsymbol{r}}_d-{\boldsymbol{r}}\right)
\end{eqnarray}
where $K_3,\dots,K_6$ can be chosen such that the roots of the characteristic equation 
\begin{eqnarray}
    \lambda^4+\lambda^3K_3+\lambda^2K_4+\lambda K_5+K_6 =0,
\end{eqnarray}
are located in the open left half complex plane. Hence, $\boldsymbol{r}$ converges to $\boldsymbol{r}_d$.

In order to find the relation between $\boldsymbol{s}$ and $\Tilde{\boldsymbol{u}}$, we need to find $\ddddot{\boldsymbol{r}}$ by differentiating twice with respect to time from $\ddot{\boldsymbol{r}}$. From~\eqref{non-linear state space}, we have
\begin{eqnarray}
    \ddot{\boldsymbol{r}} = \frac{p}{m}\hat{\boldsymbol{k}}_b-g\hat{\boldsymbol{k}}
\end{eqnarray}
By differentiating the above expression,
\begin{eqnarray}
    \dddot{r} = \frac{\Dot{p}}{m}\hat{\boldsymbol{k}}_b + \frac{p}{m} \boldsymbol{\omega}\times\hat{\boldsymbol{k}}_b
\end{eqnarray}
\begin{eqnarray}\label{ddddotr}
    \ddddot{r} = \frac{1}{m}(\boldsymbol{O}_1\boldsymbol{\Theta}+\boldsymbol{O}_2),
\end{eqnarray}
where $\boldsymbol{\Theta}=\begin{bmatrix}\ddot{p}&\ddot{\phi}&\ddot{\theta}&\ddot{\psi}\end{bmatrix}^T$ and
\begin{eqnarray}
    \boldsymbol{O}_1 = \begin{bmatrix}
     \hat{\boldsymbol{k}}_b& -p\hat{\boldsymbol{j}}_b& p(\hat{\boldsymbol{j}}_2 \times \hat{\boldsymbol{k}}_b)& p(\hat{\boldsymbol{k}}_1)\times \hat{\boldsymbol{k}}_b 
    \end{bmatrix}
\end{eqnarray}
\begin{eqnarray}
    \boldsymbol{O}_2 = p\Tilde{\boldsymbol{B}_2}\times \hat{\boldsymbol{k}}_b + \boldsymbol{\omega}\times(\boldsymbol{\omega\times p \hat{\boldsymbol{k}}_b})+2\Dot{p}\boldsymbol{\omega\times\hat{\boldsymbol{k}}_b}
\end{eqnarray}
where $\Tilde{\boldsymbol{B}_2}$ is defined in~\eqref{tilde_B2}. From~\eqref{equation: u to ddot}, $\boldsymbol{\Theta}$ can be written in the form of
\begin{eqnarray}\label{Theta}
    \boldsymbol{\Theta} = \boldsymbol{O}_3\Tilde{\boldsymbol{u}} + \boldsymbol{O}_4,
\end{eqnarray}
where
\begin{eqnarray}
    \boldsymbol{O}_3 =  
    \begin{bmatrix}
    1 & \boldsymbol{0}_{1\times3}\\
    \boldsymbol{0}_{1\times3} & \boldsymbol{J}^{-1}{\Tilde{\boldsymbol{B}}_1}^{-1} \end{bmatrix},
\end{eqnarray}
\begin{eqnarray}
    \boldsymbol{O}_4 =  
    \begin{bmatrix}
    0\\
    -{\Tilde{\boldsymbol{B}}_1}^{-1}\Tilde{\boldsymbol{B}}_2-\boldsymbol{J}^{-1}\boldsymbol{\omega}\times(\boldsymbol{J}\boldsymbol{\omega}) \end{bmatrix}.
\end{eqnarray}
Substituting \eqref{Theta} in \eqref{ddddotr}
\begin{eqnarray}\label{control input}
    \boldsymbol{s} = \frac{1}{m} \left( \boldsymbol{O}_1\boldsymbol{O}_3\Tilde{\boldsymbol{u}} + \boldsymbol{O}_1\boldsymbol{O}_4+\boldsymbol{O}_2\right)
\end{eqnarray}
%-----------------------------------------
\begin{algorithm}
\caption{Trajectory Recovery Algorithm for $i^\text{th}$ UAV }\label{alg:cap}
\begin{algorithmic}
    \State \textbf{Input} $\boldsymbol{r}_{i,0}, \mathcal{F}$ and $\omega^{\max}_{r,i}$
    \State \textbf{Output}
    $\boldsymbol{r}_{i,d}(t)$,$\boldsymbol{\omega}_{r,i}(t)$
    \State Discretize the environment in $x-y$ plane
    \State Compute $\Psi$ from \eqref{eq: L * phi=0}
    \State Compute stream lines ($\Psi$ constant curves)
    \State Find $\boldsymbol{r}_{i,d}(t)$ as a contour line corresponding to $\Psi_{i,0}$
    \While{$\omega_{r,i} < \omega^{\max}_{r,i}$}
        \State Increase quadcopter's translation speed
        \State Find control input from \eqref{control input}
        \State Compute $\boldsymbol{\omega}_{i,r}(t)$ from \eqref{eq: angular rotor speed}
    \EndWhile
\end{algorithmic}
\end{algorithm}

%-----------------------------------------
\section{Simulation Results}\label{Simulation Results}
In this section, we deploy the proposed recovery and control approach to the motion planning of a group of quadcopters. We consider a given airspace $\mathcal{C}$, in which a set of failures $\mathcal{F}=\{n_1,n_2\}$ is reported at specific positions $\boldsymbol{r}_{n_1},\boldsymbol{r}_{n_2}$. We consider a group of 10 similar quadcopters at different positions at $t_0$ (Fig~\ref{fig:initial conf}). Quadcopters' specification are listed in Table~\ref{Table}. In this scenario, all agents should modify their trajectories such that the collision avoidance and safety conditions are satisfied. To do so, we consider each failure zone as a circular cylinder of radius 2 and centered at $(x_{n_1},y_{n_1}),(x_{n_2},y_{n_2})$ along $z$-axis direction. Note that collision avoidance are guaranteed by the recovery trajectories obtained from the potential function and stream lines in Section~\ref{High-Level Planning: Recovery Trajectory Planning}. 

\begin{table}[h!]
\centering
 \begin{tabular}{|c |c| c| c|} 
 \hline
 $m$ & $g$ & $l$ & $I_x$ \\  
 \hline
 0.468 & 9.81 & 0.225 & $4.856\times 10^{-3}$\\
 \hline \hline
 $I_y$ & $I_z$ & $b$ & $k$\\ 
 \hline
 $4.856\times 10^{-3}$ & $8.801\times 10^{-3}$ & $2.98\times 10^{-6}$ & $1.14\times 10^{-7}$\\
 \hline
 \end{tabular}
 \caption{Quadcopters' specification}
 \label{Table}
\end{table}

Using the proposed technique in Section~\ref{High-Level Planning: Recovery Trajectory Planning} enables to update the trajectory of each agent, based on the stream function $\Psi$ over $\mathcal{C}$. Fig~\ref{fig:stream_contour} shows the contours of $\Psi$ constants in $x-y$ plane. 

%\begin{figure}[ht]
%    \centering
%    \includegraphics[width=0.45\textwidth]{plots/2D.png}
%    \caption{Quadcopters configuration in $x-y$ plane: blue at $t_0=0 \sec$, green at $t_0=20 \sec$, yellow at $t_0=40 \sec$ and magenta at $t_0=60 \sec$. Red diamonds show the position of the failures.}
%    \label{fig:stream_contour}
%\end{figure}

\begin{figure}[ht]
    \centering
    \includegraphics[width=0.45\textwidth]{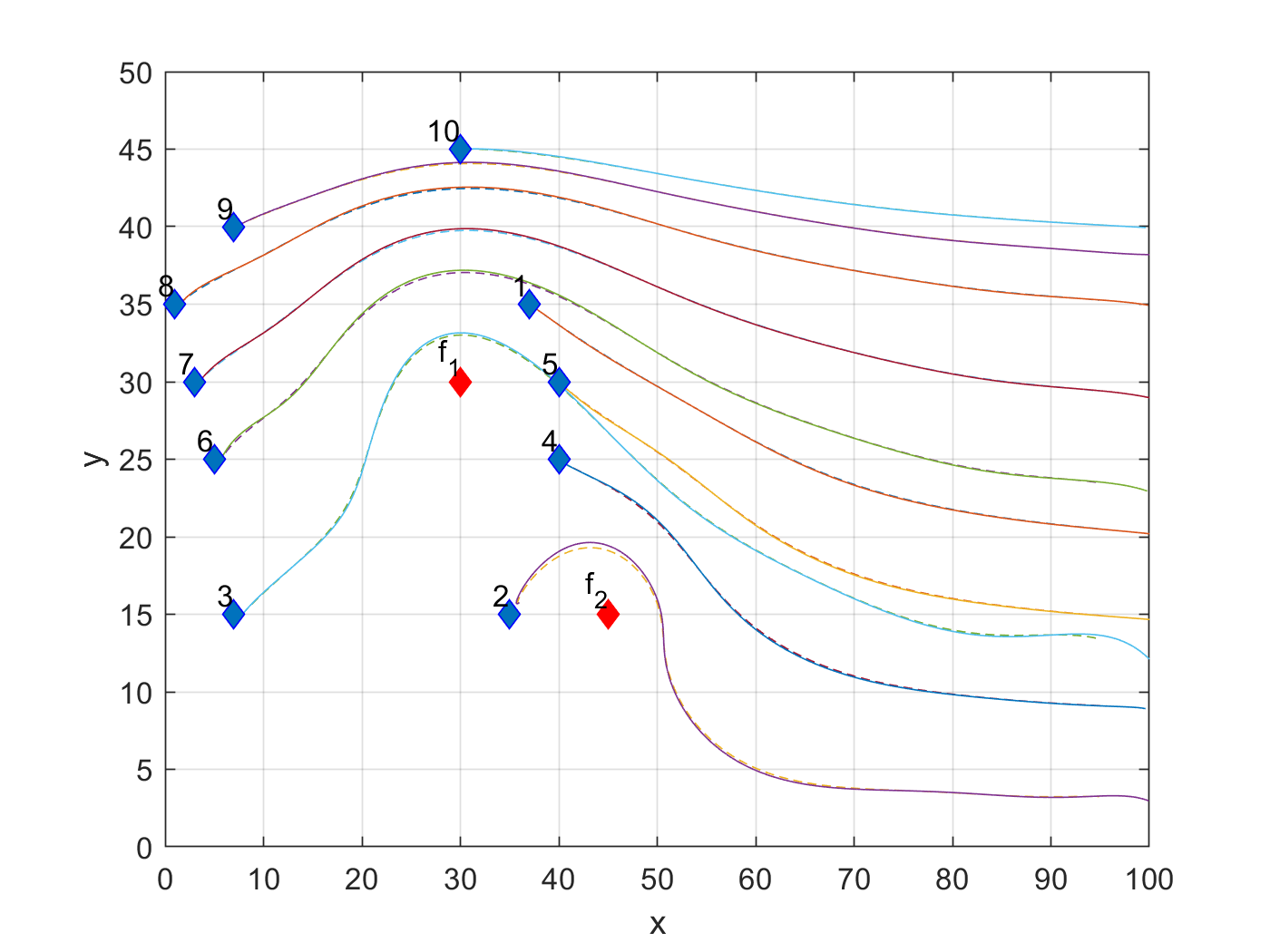}
    \caption{Solid lines show the desired trajectory $\boldsymbol{r}_{i,d}$, and dashed lines show the actual trajectory of each quadcopter.}
    \label{fig:trajectory}
\end{figure}

In the next step, desired trajectory $\boldsymbol{r}_{i,d}(t)$ is assigned to each quadcopter based on the initial position and~\eqref{eq: psi initial}. As shown in Fig~\ref{fig:stream_contour}, desired trajectories are smooth functions. We use the curve-fitting toolbox in MATLAB to approximate the desired trajectory as a polynomial function in $x-y$ plane, and consequently, we figure out the time derivatives of corresponding desired trajectories. Fig~\ref{fig:trajectory} shows the desired trajectories and the actual trajectories of each quadcopter by using the control input proposed in Section~\ref{Mathematical Modeling of Quadcopters and Trajectory Tracking Control}. We choose $K_1=1,K_2=1,K_3= 14, K_4=71,K_5=154$ and $K_6=120$ as control parameters.

\begin{figure}[ht]
    \centering
    \includegraphics[width=0.45\textwidth]{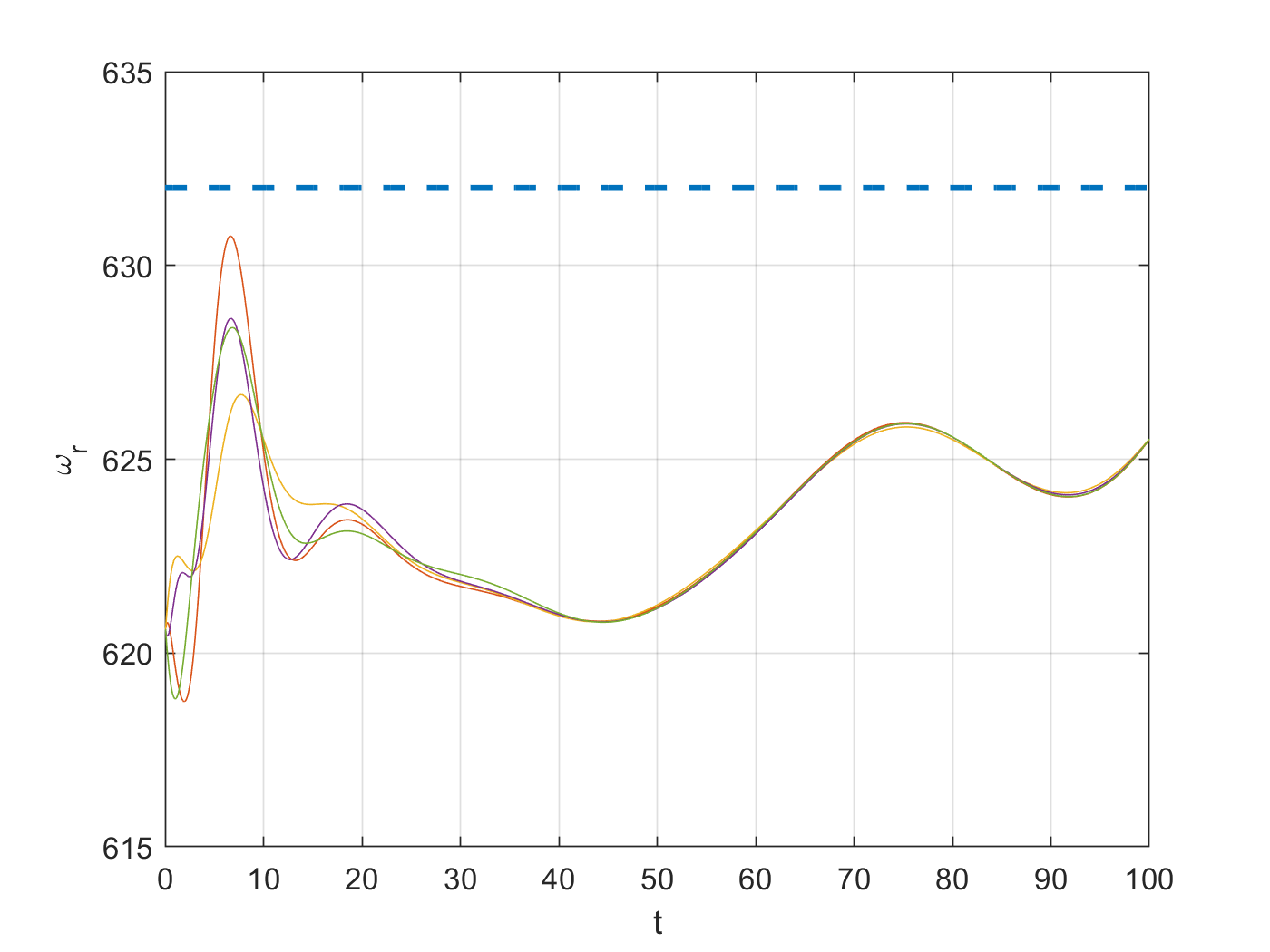}
    \caption{Angular speeds of quadcopter 2. Dashed line shows $\omega_r^{\max}$}
    \vspace{-0.5cm}
    \label{fig:rotor_angular_speed}
\end{figure}

In order to satisfy the safety condition~\eqref{mainsafety} and keep the angular speed of rotors in the safe performance limit, translation speed of each agent $i$ can be changed along a desired trajectory of $\boldsymbol{r}_{d,i}(t)$. Thus, the finite horizon optimal problem can be solved numerically to find the optimal speed for each quadcopter such that $\omega_r{_{i,j}} < \omega^{\max}_r$ for $j=1,\dots,4$. Fig~\ref{fig:rotor_angular_speed} shows the angular speeds of quadcopter which is upper-bounded by $\omega^{\max}_r$.

\section{Conclusion}\label{Conclusion}
We developed a new physics-based method for fault-resilient multi-agent coordination in the presence of unpredictable agent failure. Without loss of generality, we assumed that agents represent quadcopters that are modeled by $14$-th order nonlinear dynamics. By classifying quadcopters as healthy and failed agents, coordinating in a shared motion space, we defined the safety recovery paths of the healthy quadcopters as streamlines in an ideal fluid flow wrapping failed quadcopters. To assure quadcopter coordination safety is recovered as quickly as possible, desired trajectories of cooperative quadcopters were determined by maximization of sliding speed along the recovery streamlines such that rotor speeds of all quadcopters do not exceed a certain upper limit at all times. We also show that every healthy quadcopter can stably track the desired recovery trajectory by applying the input-output feedback linearization control.

\section{Acknowledgement}
This work has been supported by the National Science
Foundation under Award Nos. 2133690 and 1914581.

\bibliographystyle{IEEEtran}
\bibliography{ref}
\end{document}